\documentclass[letterpaper, 10 pt, conference]{ieeeconf}

\IEEEoverridecommandlockouts                              
 \pdfoutput=1
\usepackage{amsmath, amssymb}
\usepackage{amsfonts}
\usepackage{graphicx}
\usepackage{enumerate}
\usepackage{ifthen}
\usepackage{multicol}
\usepackage{float}
\usepackage{cite}
\usepackage{fancyhdr}
\usepackage[usenames, dvipsnames]{color}
\graphicspath{{figures/}}
\usepackage{mathtools}
\usepackage{multicol}
\usepackage{xcolor}

 \usepackage{cancel} 

\newtheorem{theorem}{Theorem}
\newtheorem{lemma}[theorem]{Lemma}
\newtheorem{proposition}[theorem]{Proposition}
\newtheorem{corollary}[theorem]{Corollary}
\newtheorem{rem}{Remark}

\newboolean{showcomments}
\setboolean{showcomments}{false}

\newcommand{\todo}[1]{  \ifthenelse{\boolean{showcomments}}
{\textcolor{ForestGreen}{TO DO:  #1}}{}}
\newcommand{\richard}[1]{\ifthenelse{\boolean{showcomments}}
{\textcolor{Orange}{(Richard says: #1)}}{}}
\newcommand{\alain}[1]{\ifthenelse{\boolean{showcomments}}
{\textcolor{Blue}{(Alain says: #1)}}{}}
\newcommand{\jonas}[1]{\ifthenelse{\boolean{showcomments}}
{\textcolor{Blue}{(Jonas says: #1)}}{}}
\newcommand{\kristian}[1]{\ifthenelse{\boolean{showcomments}}
{\textcolor{Blue}{(Kristian says: #1)}}{}}
\newcommand{\emma}[1]{\ifthenelse{\boolean{showcomments}}
{\textcolor{VioletRed}{(Emma says: #1)}}{}}
\newcommand{\ifneeded}[1]{\ifthenelse{\boolean{showcomments}}
{\textcolor{Gray}{#1}}{}}

\newboolean{showedit}
\setboolean{showedit}{true}
\newcommand{\edit}[1]{\ifthenelse{\boolean{showedit}}
{\textcolor{Blue}{#1}}{}}
\newcommand{\draft}[1]{\ifthenelse{\boolean{showedit}}
{\textcolor{gray}{#1}}{}}

\newcommand{\ddt}{\frac{\mathrm{d}}{\mathrm{d}t}}

\usepackage{array}
\newcolumntype{L}[1]{>{\raggedright\let\newline\\\arraybackslash\hspace{0pt}}m{#1}}
\newcolumntype{C}[1]{>{\centering\let\newline\\\arraybackslash\hspace{0pt}}m{#1}}
\newcolumntype{R}[1]{>{\raggedleft\let\newline\\\arraybackslash\hspace{0pt}}m{#1}}

\usepackage{caption}
\usepackage{subcaption}
\usepackage{tikz}
 \usetikzlibrary{plotmarks}

\usepackage{tikz}
\usepackage{pgfplots}
\pgfplotsset{compat=newest}
\usetikzlibrary{patterns}
\usetikzlibrary{decorations.text}
\usepgfplotslibrary{fillbetween}

\usepackage{flushend}
\renewcommand{\baselinestretch}{.97}

\providecommand{\figref}{}
\renewcommand{\figref}[1]{Fig.~\ref{fig:#1}}
\providecommand{\secref}{}
\renewcommand{\secref}[1]{Sec.~\ref{sec:#1}}

\usepackage{lipsum}
\usepackage{mathtools}

\title{\LARGE \bf Limitations of time-delayed case isolation in heterogeneous SIR models}

\author{ {Jonas Hansson, Alain Govaert, Richard Pates, Emma Tegling, Kristian Soltesz} 
 \thanks{The authors are with the Department of Automatic Control,
        Lund University, Lund, Sweden. Email: \{{\tt\small{jonas.hansson,  alain.govaert, richard.pates, emma.tegling, kristian.soltesz}\}@control.lth.se}}\\
        \thanks{This work was partially funded by the Swedish Research Council (grants 2017-04989, 2019-00691) and by the Wallenberg AI, Autonomous Systems and Software Program (WASP) funded by the Knut and Alice Wallenberg Foundation. The authors are members of the ELLIIT Strategic Research Area at Lund University.}
}

\begin{document}
\maketitle
\begin{abstract}
Case isolation, that is, detection and isolation of infected individuals in order to prevent spread, is a strategy to curb infectious disease epidemics. Here, we study the efficiency of a case isolation strategy subject to time delays in terms of its ability to stabilize the epidemic spread in heterogeneous contact networks. For an SIR epidemic model, we characterize the stability boundary analytically and show how it depends on the time delay between infection and isolation as well as the heterogeneity of the inter-individual contact network, quantified by the variance in contact rates. We show that network heterogeneity accounts for a restricting correction factor to previously derived stability results for homogeneous SIR models (with uniform contact rates), which are therefore too optimistic on the relevant time scales.
We illustrate the results and the underlying mechanisms through insightful numerical examples. 
\end{abstract}

\section{Introduction}
\subsection{Background}
Early in an epidemic of a previously unknown infectious disease, transmitted through an unknown pathogen, existing pharmaceutical treatments may not be applicable or efficient. This was for example the case with COVID-19, caused by the pathogen SARS-CoV-2. In early 2020, countries across the world experienced exponential increases in detected SARS-CoV cases, hospitalizations, and ICU admittance due to COVID-19, and ultimately related deaths. While statistics are still debated one year later~\cite{oliver21}, it stood clear early on that pharmaceutical treatment and clinical care would not alone suffice to curb the epidemic outbreak.

The exponential increase of cases early in the epidemic prompts the implementation of schemes to reduce the number of probable transmission events. Such schemes are often referred to as non-pharmaceutical interventions (NPIs for short), and can be partitioned into two categories:
\begin{itemize}
    \item Recommendations or legislation aimed at decreasing inter-individual contact rates;
    \item Schemes for case isolation through testing, and possibly contact tracing.
\end{itemize}

The former includes, for example, discouraging unnecessary in-person interaction, banning large gatherings, closing schools, or imposing societal lockdowns. The effectiveness of such interventions is still poorly understood. In particular, it can be practically challenging or impossible to estimate the effect of individual NPIs, particularly in the early transient phase of an epidemic, as demonstrated in \cite{soltesz20,gustafsson21}.

In light of the above, case isolation schemes have gained attention. In contrast to NPIs aimed at
reducing contact rates in general, case isolation schemes target contacts involving infectious individuals. If successfully implemented, they could be the key to an open society where lockdowns are replaced by regular testing. Under such schemes, a positive test result leads to case isolation. The effectiveness of the scheme can be increased by complementing systematic testing with contact tracing.


It is intuitively evident that a case isolation scheme needs to identify and isolate positive cases in a timely manner in order to be effective. In~\cite{pates21}, a comprehensive analysis was done, from a control-theoretic standpoint, of the performance and stability of case isolation schemes. We refer to references therein for relevant works in the epidemiology literature.  There, as in our current work, one central question is that of time delays: how long can delays between infection and isolation be if we are to prevent epidemic outbreaks?  

In this paper we investigate how the structure of social networks affect the effectiveness of case isolation schemes. This is done by deriving an explicit stability 
condition, which relates effectiveness to network heterogeneity. We use a simple epidemiological model; an SIR model for the early phases of an epidemic. This modeling choice is
motivated by its sufficiency for the purpose of our analysis, and the lack of data to support high-complexity modeling.

\subsection{Preliminaries}
We consider a closed population. This is an optimistic approximation, since there will be no imported cases within the model. The population is partitioned into susceptibles, infectious and removed (recovered and permanently immune or deceased) proportions. Here, we will not distinguish between being infected and being infectious. However, the results we are about to present can easily be extended to cater for this distinction.

If the proportion of infectious goes to zero, the strategy has been successful.
The question we consider is how large the infectious proportion can grow, relative to the susceptible, in order for a case isolation strategy to remain effective.

To answer this, we consider the setting early in an epidemic, where the proportion of susceptibles is much larger than the infectious and removed, respectively. This allows us to disregard herd immunity effects. Another way to put it is that in an interaction between two individuals where one is infectious, the other is susceptible (with probability one).

Furthermore, we realistically assume the considered time window for our model to be short compared to the typical duration of immunity. This means that we neglect any flow from the removed to the susceptible sub-population. As with distinguishing between infected and infectious individuals, introducing such reflux dynamics into the model we propose is straightforward, should the considered disease differ in this aspect from SARS-CoV-2.

\emma{Add sth on limitations of the model and the linearization here. We analyze a worst-case scenario and want to stabilize the disease-free equilibrium. With this model we get conservative stability conditions. 
}
We will focus on criteria for stabilization of the early epidemic trajectory. We will refer to a case isolation scheme as stabilizing if it eventually empties the infectious sub-population. In the linearized model, this is done without decreasing the pool of susceptibles (in other words without the help of  the herd immunity effect~\cite{topley23}), and in this sense, the stability criteria we derive are conservative. 


Note that this stability condition does not specify \emph{performance} in the sense that it 
may allow for a large infectious sub-population in the transient during which the infectious sub-population grows (due to the time delay before isolation) and empties. 
As such, this stability condition constitutes a bare minimum. Any practically feasible strategy would need to fulfill it with some performance margin, as further elaborated in \cite{pates21}. 

Finally, we need to formalize our case isolation scheme. Here we will utilize a simple yet versatile model that quantifies the proportion of those infected on a given day that is subsequently isolated by the scheme $T_\text{delay}$ days later. 
This can be readily reparameterized into, for example, the frequency with which individuals are tested, compliance to the scheme, logistic delays and specificity of the employed test for infection.

\section{Case isolation in heterogeneous populations}\label{sec:heterogeneous}
Assuming that the population is large, so that quantization effects become negligible, 
the trajectory of the epidemiological state subject to case isolation can be modeled as in \cite{pates21}:
\begin{equation}\label{eq:homogeneous SIR}
\ddt
\begin{bmatrix}
S\\I\\R
\end{bmatrix}
=
\begin{bmatrix}
-1\\1\\0
\end{bmatrix}
\beta S (I-Q) + 
\begin{bmatrix}
0\\-1\\1
\end{bmatrix}
\gamma I,\quad
S+I+R=1.
\end{equation}
This model, which we will term \emph{homogeneous} because of the uniform contact rates across the entire population, 
is a slight variant of the traditional SIR model of \cite{KM27}. Here, as usual, $S$, $I$ and $R$ denote the proportions of the population that are susceptible, infectious, and removed and, as throughout, we have dropped the time-dependence in the variables to simplify notation when possible. The transitions between the states are governed by two rates that model the effect of disease spread and recovery, with mixing and recovery parameters $\beta$ and $\gamma$, respectively. Note in particular that the mixing rate has been adjusted to account for the effect of isolating infectious individuals.
More specifically, 
\begin{equation}\label{eq:feedback}
Q(t)=\alpha{}e^{-\gamma{}T_{\text{delay}}}I(t-T_{\text{delay}})
\end{equation}
denotes the proportion of the population that is both infectious and isolated, subject to a time delay $T_{\text{delay}}$ after infection. The rate that describes the spread of the disease has been modified so that the spread is only driven by interactions between the remaining infectious population and the susceptible population (the $\beta{}S(I-Q)$ term). 

This captures to a first approximation the two most important features of a case isolation scheme, namely, how quickly infectious individuals are identified and isolated ($T_{\text{delay}}$), and what proportion of cases are found ($\alpha$).
However, the homogeneously interacting population is clearly a simplification.
\color{black}
In reality, the interactions that may lead to infection are much more complex and hard to model. For example, you interact more frequently with people at your workplace than with people from a remote country.
There is also a time-varying aspect. For instance, if your workplace issues a work-from-home guideline to reduce disease transmission, your rate of interaction with your coworkers will typically drop.

Numerous works have been dedicated to modeling the associated contact network dynamics, with \cite{may06,salathe10} and survey~\cite{nowzari2016analysis} constituting representative examples. The validity of such network models is hard to verify, and the time-varying aspects further increase the uncertainty surrounding their accuracy. We therefore delimit ourselves to a simple but important question: how does the introduction of interaction heterogeneity alter the requirements on our case isolation scheme \eqref{eq:feedback}, expressed in terms of the isolation proportion parameter $\alpha$ and associated time delay $T_\text{delay}$?

\subsection{Heterogeneous population model}

To account for contact rate heterogeneity, early models of infectious diseases (particularly STDs)~ \cite{may1987commentary,may2001infectious} incorporated a distribution of contact rates $N_k/N$ defined as the proportion of the population of size $N$, who on average have $k$ contacts per time unit (and in all other regards are homogeneous). It is assumed that a contact between a susceptible and infectious individual transmits an infection at a rate $\rho$.
Using the formulation of \cite{anderson1986preliminary}, the disease dynamics for each partition $k=0,1,\dots,n$ becomes
\begin{equation}\label{eq:het model}
    \ddt
    \begin{bmatrix}
    X_k\\
    Y_k
    \end{bmatrix}=
    \begin{bmatrix}
    -k\lambda&0\\
    k\lambda&-\gamma
    \end{bmatrix}
    \begin{bmatrix}
    X_k\\
    Y_k
    \end{bmatrix},
\end{equation}
where $X_k$ and $Y_k$ denote, respectively, the number of susceptible and infectious individuals in each partition. This generalizes the standard SIR model for the homogeneous through the variable
\begin{equation}\label{eq: lambda}
  \lambda=\rho\frac{\sum_{k=1}^n kY_k}{\sum_{k=1}^n kN_k}, 
\end{equation}
which is the rate at which infection is acquired from any one randomly chosen contact per time unit---now more likely to come from the partition with higher contact rates.

\begin{rem}[Retrieving the homogeneous SIR model] 
When all individuals have $k$ contacts per time unit, the rate of acquiring an infection from a randomly chosen contact in~\eqref{eq: lambda} simplifies to $\lambda=\rho Y/N=\rho I$. The disease dynamics~\eqref{eq:het model} recover the traditional SIR model of \cite{KM27}  with mixing parameter $\beta=\rho k$ when divided by $N$.
That is,
\begin{equation}
    \ddt
    \begin{bmatrix}
    S\\
    I
    \end{bmatrix}=
    \begin{bmatrix}
    -\beta I&0\\
    \beta I&-\gamma
    \end{bmatrix}
    \begin{bmatrix}
    S\\
    I
    \end{bmatrix},\quad R=1-I-S,
\end{equation}
where $S=X/N$, $I=Y/N$ (this is equivalent to \eqref{eq:homogeneous SIR} with $Q\equiv{}0$).
\end{rem}


\subsection{Case isolation in a heterogeneous population}\label{sec: het model}
In arguably the simplest generalization of the case isolation scheme~\eqref{eq:feedback} to the heterogeneous population model~\eqref{eq:het model}, the number of individuals that are both isolated and infectious individuals in partition $k$ is equal to some common proportion $0\leq\alpha\leq1$ of the infectious individuals $Y_k$,  $T_{\text{delay}}$ days in the past.
This corresponds to a uniform identification and isolation scheme, independent of an individuals' number of contacts. Furthermore, it is assumed the partition size $N_k$ is constant: isolated cases do not leave partition $k$, but are simply excluded from transmitting the disease. 
The number of infectious individuals in partition $k$ then satisfies
\begin{equation}\label{eq: hetereogeneous case isolation}
    Y_k(t)-Q_k(t), \quad Q_k(t)=\alpha e^{-\gamma T_{\text{delay}}}Y_k(t-T_{\text{delay}}),
\end{equation}
and the rate at which an infection is acquired from any one randomly chosen contact becomes
\begin{equation}\label{eq: lambda adjusted}
  \xi(t)=\rho\sum_i i(Y_i(t)-Q_i(t))/\sum_k kN_k.
 \end{equation}
As a consequence of \eqref{eq: lambda adjusted}, the disease dynamics in each partition are described by the $n$ delayed differential equations: 
 \begin{equation}\label{eq: delayed collective dynamics}
\ddt \!\! \begin{bmatrix}
Y_1\\
Y_2\\
\vdots\\
Y_n
\end{bmatrix}\!\! =\! \frac{\rho}{\sum\limits_kN_k} \!\! \begin{bmatrix}
X_1\cdot1\\
X_2\cdot 2\\
\vdots\\
X_n\cdot n
\end{bmatrix}\!\!\!
 \begin{bmatrix}
 1\\2\\\vdots\\n
 \end{bmatrix}^{\! \top}
\!\! \begin{bmatrix}
Y_1-Q_1\\
Y_2-Q_2\\
\vdots\\
Y_n-Q_n
\end{bmatrix} 
\!-\!\gamma \mathbf{I}_n \!\!\begin{bmatrix}
Y_1\\
Y_2\\
\vdots\\
Y_n
\end{bmatrix}\!\!.
\end{equation}
About the all-susceptible equilibrium ($X_k=N_k$ for all $k$), the dynamics of the \emph{total} number of infectious individuals in the heterogeneous population, $Y=\sum_kY_k$, reads as
\begin{equation}\label{eq: Y dynamics}
\begin{aligned}
\ddt Y
&=\sum_k kN_k\xi-\gamma Y,
\end{aligned}
\end{equation}
which can be written as a delayed differential equation in terms of $\lambda$. 
In turn, by differentiating \eqref{eq: lambda} one obtains
\begin{equation}\label{eq: lambda dynamics}
         \frac{d}{dt}\lambda
         =\frac{\rho}{\sum_kkN_k}\sum_i \left(i^2N_i\xi\right)-\lambda\gamma Y.
\end{equation}
All dynamics are described in terms of the number of individuals. By dividing by $N$, 
the population level disease dynamics with delayed case isolation are described by:
\begin{multline}\label{eq: two dimensional system}
    \frac{d}{dt}\begin{bmatrix}
    I\\\lambda\end{bmatrix}=
    \underbrace{\begin{bmatrix}
    -\gamma&\mu\\
    0&\rho\left(\frac{\sigma^2}{\mu}+\mu\right)-\gamma
    \end{bmatrix}}_{\mathbf{A}_0}
    \begin{bmatrix}
    I\\\lambda
    \end{bmatrix}\\+\underbrace{\begin{bmatrix}
    0 &-\mu \alpha e^{-\gamma T_{\text{delay}}}\\
    0& -\rho\left(\frac{\sigma^2}{\mu}+\mu\right)\alpha e^{-\gamma T_{\text{delay}}} \end{bmatrix}}_{\mathbf{A}_1}\begin{bmatrix}
    I(t-T_{\text{delay}})\\
    \lambda(t-T_{\text{delay}})
    \end{bmatrix},
\end{multline}
where $\mu$ and $\sigma$ represent the mean and standard deviation of the number of contacts per time unit.
For a given delay~$T_{\text{delay}}$, the system~\eqref{eq: two dimensional system} is asymptotically stable if and only if all of the roots of 
\begin{equation}\label{eq: characteristic equation}
    \mathrm{det}\underbrace{\begin{bmatrix}s+\gamma&-\mu(1-\alpha e^{-(\gamma+s) T_{\text{delay}}})\\[6pt]
    0& s-\rho\left(\frac{\sigma^2}{\mu}+\mu\right)(1-\alpha e^{-(\gamma+s)T_{\text{delay}}})+\gamma\end{bmatrix}}_{\mathbf{M}=s\mathbf{I}_2-\mathbf{A}_0-\mathbf{A}_1e^{-sT_{\text{delay}}}}=0,
\end{equation}
are in the open left-half complex plane~\cite{sipahi2011stability}.
Because $\mathbf{M}$ is upper triangular, the roots of the characteristic equation~\eqref{eq: characteristic equation} are $-\gamma$ and those of
\begin{equation}\label{eq: fs}
    f(s)=s-\rho\left(\frac{\sigma^2}{\mu}+\mu\right)\left(1-\alpha \mathrm{exp}(-(\gamma+s)T_{\text{delay}})\right)+\gamma=0.
\end{equation}
Here, $f(s)$ is equivalent to the characteristic equation of the homogeneous case~\eqref{eq:homogeneous SIR} with a scaled mixing parameter $\beta$ linearized about $(S,I,R,Q)=(1,0,0,0)$. This,
together with our derivation above, 
leads to the following lemma.

\begin{lemma}\label{lem: equivalence}
About the all-susceptible equilibrium {$S_k=N_k$} for all $k=1,\dots,n$, the heterogeneous SIR model with delayed case isolation~\eqref{eq: hetereogeneous case isolation}--\eqref{eq: delayed collective dynamics}, and its population level disease dynamics~\eqref{eq: two dimensional system}, are asymptotically stable if and only if the linearization of the homogeneous model~\eqref{eq:homogeneous SIR}--\eqref{eq:feedback} about $(S,I,R,Q)=(1,0,0,0)$ is asymptotically stable with mixing parameter $\beta=\rho\mu(c_v^2+1)$. 
\end{lemma}

We will provide an explicit stability condition shortly, after reiterating the corresponding condition for the homogeneous case. 

\subsection{Comparing the homogeneous and heterogeneous population models}
In the absence of case isolation, the equivalence of the stability condition between the heterogeneous and homogeneous SIR model, early in the epidemic, was first established in~\cite{may1987commentary}. 
We now know the stability equivalence is also valid under the case isolation scheme with uniform identification and isolation of infectious individuals\footnote{Having a common $\alpha$ in each partition is not necessary for the equivalence relation between the stability condition of the heterogeneous and homogeneous model to hold. For example, when $\alpha$ scales proportional to the number of contacts $k$ in the partition ($\alpha_k=\alpha k/n$), the stability condition of the all susceptible equilibrium is equivalent to that of a common $\alpha$ that is scaled by $\langle k^3\rangle/(n\sigma^2\mu^2)$, where $\langle k^3\rangle$ is the third raw moment of the contact distribution.}. 

To compare the stability conditions of the heterogeneous and homogeneous model in terms of epidemiological variables, we briefly discuss the most common epidemiological measure (and source of confusion) for the spread of the disease: the reproduction number.

The reproduction number $\mathcal{R}$ describes the expected number of secondary infections caused by one primary infection. It holds that
\begin{equation*}
\mathcal{R}\propto\underbrace{\left(\frac{\text{infection}}{{\text{contact}}}\right)}_\rho\cdot\underbrace{\left(\frac{\text{contact}}{{\text{time}}}\right)}_{c}\cdot\underbrace{\left(\frac{\text{time}}{{\text{infection}}}\right)}_{d},
\label{eq:R}
\end{equation*}
where
\begin{itemize}
    \item $\rho\in[0,1]$: transmission rate given contact between a susceptible and infectious individual;
    \item $c\in\mathbb{R}_+$: average contact rate between susceptible and infectious individuals;
    \item $d\in\mathbb{R}_+$: average duration of infectiousness.
\end{itemize}
For the homogeneous population model \eqref{eq:homogeneous SIR} we have that $\beta=\rho c$ and $d=1/\gamma$. It is also evident from \eqref{eq:R} that the reproduction number reveals how \emph{much} an epidemic grows, and not how \emph{fast}. In order to quantify the latter, the serial interval between infections is needed in addition. Yet, the reproduction number is more commonly used than the corresponding growth rate $r=\beta-\gamma$ among infectious disease epidemiologists, which is why we have chosen to use the former in our parametrization of fundamental limitations.

The basic reproduction number $\mathcal{R}_0$ is the reproduction in absence of a \emph{considered} intervention. There is a common misconception that for a particular pathogen $\mathcal{R}_0$ is a universal constant (that can be looked up in the literature). Instead it depends on, among other time-varying parameters, the virulence of the pathogen and the societal structure under consideration. For instance, a particular virus would typically result in different $\mathcal{R}_0$'s in two countries, or in a city versus a village.

If an intervention is enacted, it is instead common to talk about the resulting effective reproduction number $\mathcal{R}_e$ (sometimes referred to as the time-varying reproduction number $\mathcal{R}_t$), and it really only makes sense to consider $\mathcal{R}_0$ in relation to a particular intervention: $\mathcal{R}=\mathcal{R}_0$ in absence of the intervention; $\mathcal{R}=\mathcal{R}_e$ if the intervention is enacted. 

\begin{figure}
\centering
\hspace{-1cm}
\includegraphics[scale=0.93]{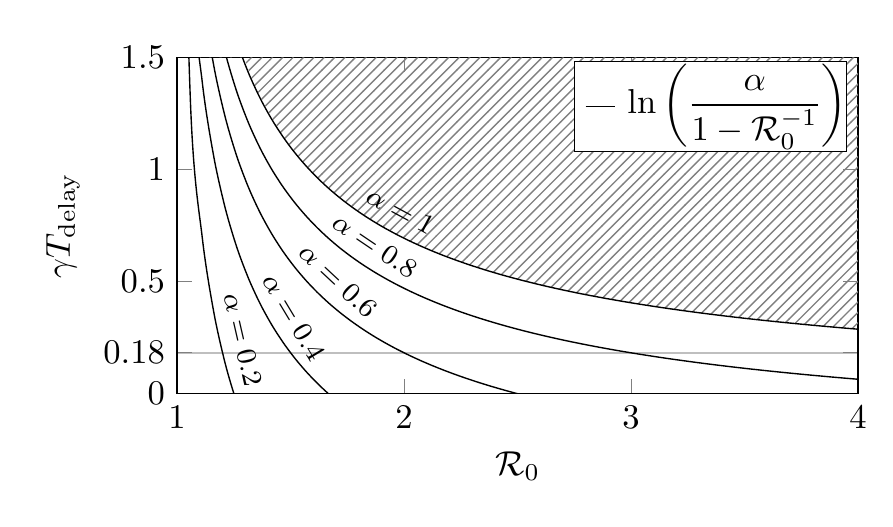}
\caption{Illustration of the stability boundary for the model \eqref{eq:homogeneous SIR}--\eqref{eq:feedback}. The model is stable if and only if $\left(\mathcal{R}_0,\gamma{}T_{\text{delay}}\right)$ lies below the corresponding $\alpha$ curve. That is, at least a proportion $\alpha$ of persons becoming infectious on a particular day need to be isolated $T_\text{delay}$ days later, given a natural recovery rate is $\gamma$~days$^{-1}$.}
\label{fig:R0tradeoff}
\end{figure}

It was shown in \cite{pates21} that for the case isolation scheme \eqref{eq:feedback} applied to the model \eqref{eq:homogeneous SIR} the relation between $\mathcal{R}_0$ and $\mathcal{R}_e$ is given by
\[\mathcal{R}_e=\mathcal{R}_0\left(1-\alpha{}e^{-\gamma{}T_{\text{delay}}}\right),\quad \mathcal{R}_0=\dfrac{\beta}{\gamma},\]
which is less than 1 if and only if~\cite[Theorem 2]{pates21}
\begin{equation}
\gamma{}T_{\text{delay}}<\ln\left(\frac{\alpha}{1-\mathcal{R}_0^{-1}}\right). \label{eq:fundamental limit}
\end{equation}
The specific trade-off between parameters and delay implied by \eqref{eq:fundamental limit} is shown in \figref{R0tradeoff}. This figure can be used to quickly assess the amount of delay that can be tolerated before instability, and hence exponential growth, occurs.

For example, with $\alpha=0.8$, $\mathcal{R}_0=3$ and $\gamma=0.1$, 
parameters chosen to be representative for SARS-CoV-2, the stability condition becomes
\begin{equation}
T_{\text{delay}}\gamma{}<0.18,\quad\Longrightarrow{}\quad{}T_{\text{delay}}<1.8\text{ days}.
\end{equation}
Clearly, short isolation times are essential when dealing with an infectious disease! We also see the importance of identifying a significant proportion of cases. When $\alpha<2/3$ (that is, the scheme detects and isolates less than two thirds of the cases) exponential growth will occur even with ${T_{\text{delay}}=0}$.
We can use the equivalence relation in Lemma~\ref{lem: equivalence} to obtain the stability condition for the heterogeneous population model.
\begin{proposition}\label{prop: het stability}
Consider the population level disease dynamics with delayed case isolation linearized around the disease-free equilibrium, which is given for the heterogeneous model~\eqref{eq: two dimensional system}. The system dynamics are stable if and only if
\begin{equation}\label{eq:stability het}
T_{\text{delay}} <\frac{1}{\gamma} \ln\left(
\frac{\alpha\beta h}{\beta h-\gamma}
\right),\quad h=(c_v^2+1), 
\end{equation}
where $c_v=\sigma/\mu$ 
is the the coefficient of variation, or relative standard deviation. 
\end{proposition}
\begin{proof}
The result follows from \cite[Theorem 2]{pates21} together with the stability equivalence relation established in Lemma~\ref{lem: equivalence}.
\end{proof}

The stability condition in Proposition~\ref{prop: het stability}
implies that contact heterogeneity ($c_v>0$) 
would \emph{increase} 
the upper bound on $T_\text{delay}$ if and only if the reproduction number of the corresponding homogeneous system ({$\beta=\rho\mu$}), {in absence of control~\eqref{eq:feedback}}, fulfils $\mathcal{R}_0<1/h$. Since in the heterogeneous model $c_v>0\Leftrightarrow h>1$, this requirement implies that heterogeneity would allow for a longer $T_\text{delay}$ 
only if the uncontrolled system is \emph{already stable} in the sense that $\mathcal{R}_0<1$. 
However, and of larger practical importance, the upper bound on the admissible $T_\text{delay}$ will \emph{decrease} when the heterogeneity ($c_v$) is increased whenever $\mathcal{R}_0>1$.

\begin{figure}
\centering
\hspace{-1cm}
\includegraphics[scale=0.93]{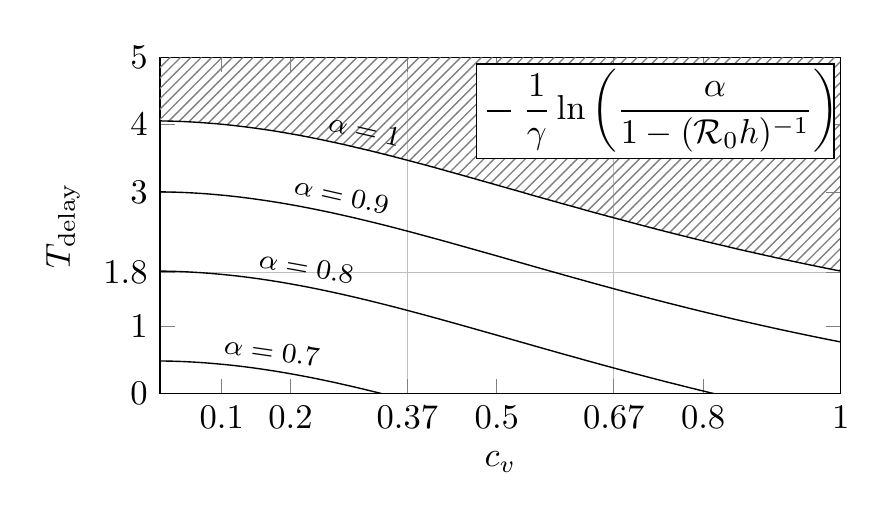}
\caption{
Illustration of the stability boundary for the model~\eqref{eq: two dimensional system}  
with $\mathcal{R}_0=3$ and $\gamma=0.1$. The model is stable if and only if $(c_v,T_{\text{delay}})$ lies below the corresponding $\alpha$ curve. When $\alpha=0.8$, the condition \eqref{eq:cvmax} reveals that the maximum coefficient of variation for which a positive isolation delay $T_\text{delay}>0$ is allowed evaluates as $c_v\approx0.82$ (where the $\alpha=0.8$ intersects the horizontal axis). The vertical lines at $c_v=0.37$ and $c_v=0.67$ represent bounds on~$c_v$ for human infectious disease transmission networks from~\cite{salathe10}.}
\label{fig:maxcv}
\end{figure}

In the latter case, one may ask how large the coefficient of variation can be to allow for a positive delay and a stable equilibrium $(I,R,Q)=(0,0,0)$ of the linearized model with the scaled parameter $\beta$. 
\begin{corollary}
For a basic reproduction number $\mathcal{R}_0>1$ of the corresponding homogeneously interacting population and $\alpha\in[0,1)$, a positive delay is allowed in the heterogeneous population if and only if the coefficient of variance of the number of contacts satisfies \begin{equation}\label{eq:cvmax} 
    c_v^2<\frac{1}{\mathcal{R}_0(1-\alpha)}-1.
\end{equation}
\end{corollary}
\vspace{2mm}
\figref{maxcv} illustrates the effect of the coefficient of variance of the number of contacts on the upper bound of the delay~\eqref{eq:stability het}, and the maximum allowed coefficient of variance $c_v$  for several~$\alpha$ that are isolated according to \eqref{eq: hetereogeneous case isolation}. 


\color{black}

\section{Numerical simulations}
In this section, we study the role of network heterogeneity in the SIR model numerically in order to illustrate the stability condition~\eqref{eq:stability het}.

\subsection{Neighbors of infectious individuals}

The phenomenon of the effective mixing parameter~$\beta$ increasing with increased heterogeneity (quantified by the coefficient of variation $c_v$) is also the answer to the question ``Why your friends have more friends than you do'', explained with mathematical insight in \cite{feld91}. 
The intuition is that individuals with many contacts are more likely to get infected, and therefore the infectious proportion of the population will comprise individuals who have more social contacts (higher degree) than those in the susceptible proportion of the population.

\begin{figure*}[!ht]
     \centering
    \begin{subfigure}[]{.3\linewidth}
        \centering
        \includegraphics[width=1.05\linewidth]{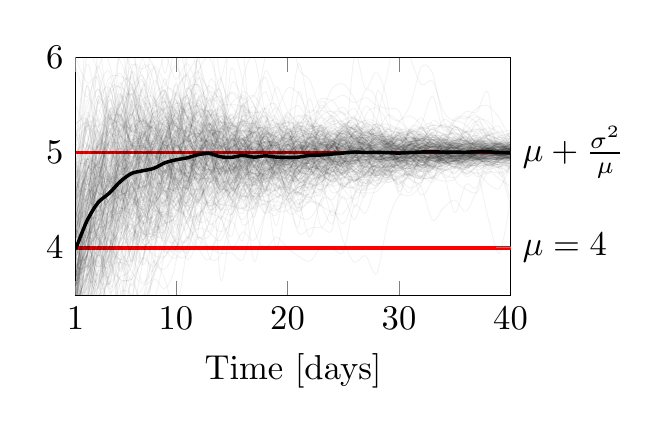}
        \caption{Configuration model with Poisson distributed degrees and uniform initialization. }
    \end{subfigure}
    \hfill
    \begin{subfigure}[]{.3\linewidth}
        \centering
        \includegraphics[width=1.05\linewidth]{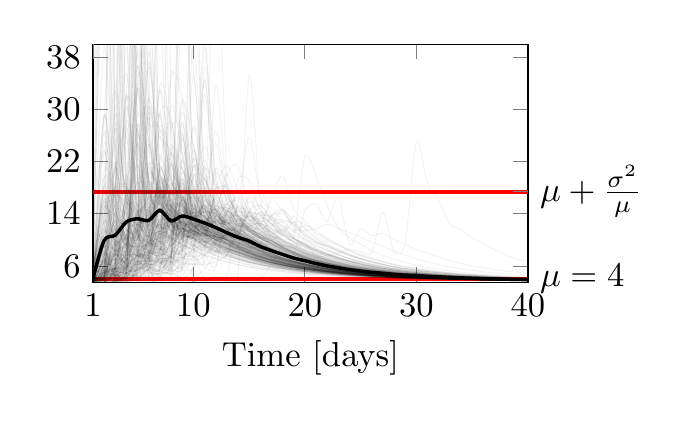}
        \caption{Barabási–Albert graph and uniform initialization}
    \end{subfigure}
    \hfill
    \begin{subfigure}[]{.3\linewidth}
        \centering
        \includegraphics[width=1.05\linewidth]{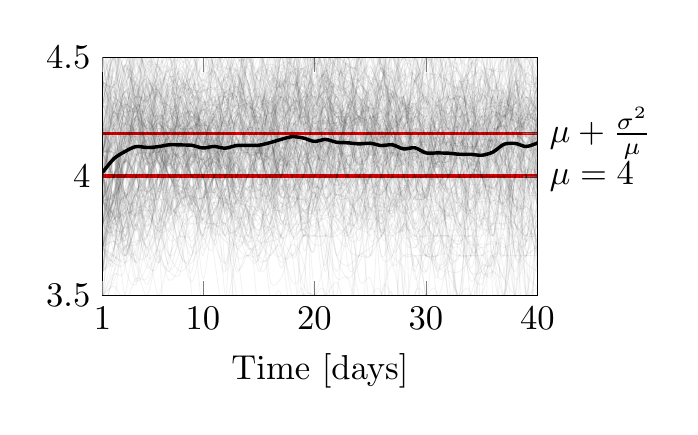}
        \caption{Watts–Strogatz graph and uniform initialization.}
    \end{subfigure}
    
    \medskip
    
    \begin{subfigure}[]{.3\linewidth}
    \centering
    \includegraphics[width=1.05\linewidth]{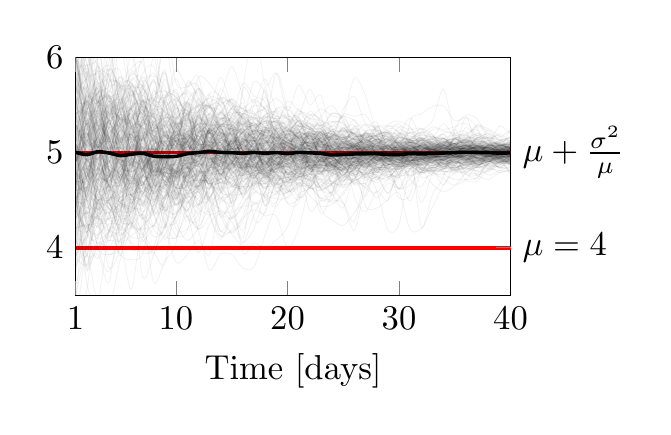}
    \caption{Configuration model with Poisson distributed degrees and degree-based initialization.} 
    \end{subfigure}
    \hfill
    \begin{subfigure}[]{.3\linewidth}
        \centering
        \includegraphics[width=1.05\linewidth]{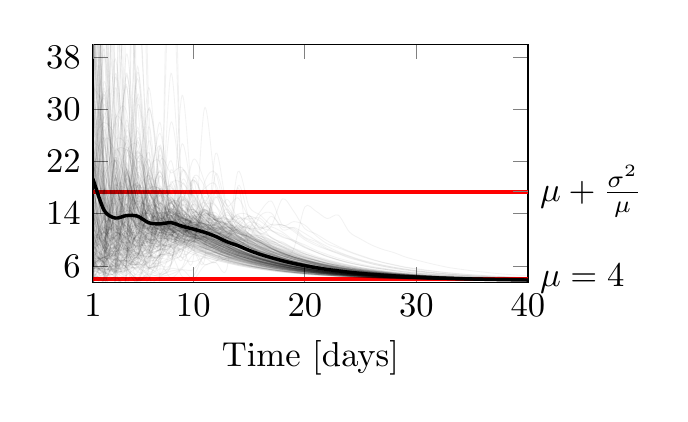}
        \caption{Barabási–Albert degree-based initialization.}
    \end{subfigure}
    \hfill
    \begin{subfigure}[]{.3\linewidth}
        \centering
        \includegraphics[width=1.05\linewidth]{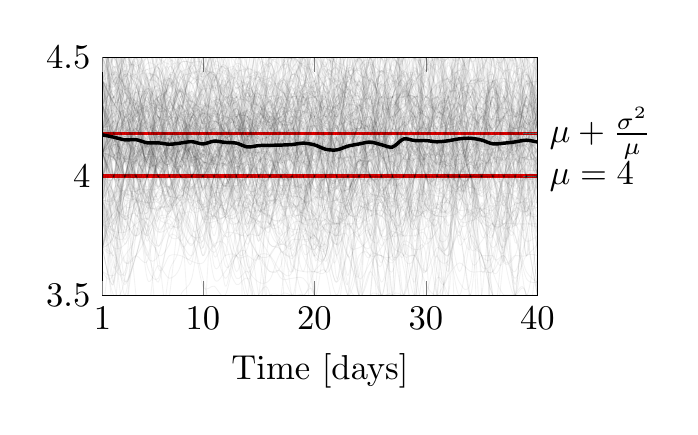}
        \caption{Watts–Strogatz graph and degree-based initialization.}  
    \end{subfigure}
     \caption{300 simulations were run for each scenario and then the average degree of the infectious population was calculated for each day. The simulation started with $10$ infectious individuals. The dark line represent the average among all simulations. The effective number of contacts $\mu + \sigma^2/\mu$ and mean of degree distribution $\mu$ are marked in each sub-figure.}
     \label{fig:rand_start_mean_neighbors}
\end{figure*}

As with friendships, the contacts considered here are mutual, translating into an undirected network graph. Thus,  the early spread of the disease is proportional to the mean number of contacts of a randomly chosen contact of an infectious individual (node), given by $\mu h$ in \eqref{eq:stability het} and illustrated in \figref{rand_start_mean_neighbors}.
Treating the network as homogeneous ($c_v=0 \Leftrightarrow h=1$) can therefore lead to an under-estimation of the reproduction number by disregarding the growth of the infectious sub-population, which is fueled by highly connected individuals, who are both more likely to acquire and spread the disease.

In \figref{rand_start_mean_neighbors} we have illustrated this effect through simulations of the early stage of 300 epidemics on different undirected random graphs representing a population. The nodes of each graph represent individuals; the edges contacts. For the example in \figref{rand_start_mean_neighbors}, we chose to generate three different random graphs based on: 1) the configuration model, 2) Barabási-Albert and 3) Watts-Strogatz; for details on these random network models, see \cite{barabasi2016}. These random graphs were chosen to represent three fundamentally different networks. 
Each graph was generated with $10^6$ nodes and an average degree $\mu=4$. 

The epidemics were each seeded by randomly assigning $10$ infectious individuals on day one. This corresponds to $I= 10^{-5} \ll S = 1-I$. The transmission probability, being the per-day probability that an edge between one infectious and one susceptible node leads to disease transmission,
was set to $\rho=0.2$.  The recovery rate from \eqref{eq:homogeneous SIR} was set to $\gamma=0.1$. Assuming independent interactions, the probability for a susceptible to become infectious the next day, given that it has $m$ infectious contacts, is thus $1 -(1-\rho)^m$.


In \figref{rand_start_mean_neighbors} (a)--(c), the results from the 300 epidemics simulations are shown. In particular, we see the average number of contacts of the infectious population throughout the early stages of the epidemic. From this figure we can clearly identify how the infectious individuals have a higher degree of contacts than the general population of the network. Furthermore, we can see that the average number of contacts of the infectious population early on tends towards a value equal to or smaller than the effective number of contacts $\mu + \sigma^2/\mu$. The same scaling defines the effective mixing parameter $\beta h$ in the stability condition for a heterogeneous population \eqref{eq:stability het}.

To test if the effective number of contacts $\mu + \sigma^2/\mu$ actually gives an upper bound for the average number of contacts of the infectious population, it is also relevant 
to simulate the epidemic starting with average degree equal to $\mu + \sigma^2/\mu$. This can be done 
by choosing the initially infectious with a probability proportional to the degree of each node. The result of this simulation is shown in \figref{rand_start_mean_neighbors} (d)--(e). Here, we note that the average degree either remains constant or decreases for the three graphs. This indicates that the scaling factor $h=\sigma^2/\mu^2+1$ is conservative and represents the infectious population early on in an epidemic.

\section{Discussion}

In \secref{heterogeneous} we generalized the case isolation scheme~\eqref{eq:feedback} of the homogeneous case isolation scheme to heterogeneous contact rate model and  characterized how the stability boundary is moved if $\mathcal{R}_0$ of a homogeneous population is used in \eqref{eq:fundamental limit}, when in fact the population is heterogeneous in the sense that the coefficient of variation of the contact network~$c_v$ is positive, while the mean degree remains unchanged.

Introducing heterogeneity in this way corresponds to multiplication of the mixing parameter $\beta$ in \eqref{eq:homogeneous SIR} with the scaling factor $h=c_v^2+1$. Since $\beta$ is directly proportional to~$\mathcal{R}_0$, the heterogeneity can be interpreted as increasing the reproduction number by the corresponding factor. Note that the same shifting applies for a population partitioned based on heterogeneity of the infectiousness parameter $\rho$. This makes it possible to apply the same modeling to account for a variability in infectiousness across variants of the considered pathogen.
Whether to multiply by the factor $h$ before using $\mathcal{R}_0$ in the analysis comes down to how $\mathcal{R}_0$ was estimated from data. If it was estimated taking the heterogeneity into account it should not be adjusted, otherwise it should. 

In the heterogeneous model \eqref{eq:het model} proposed in \cite{anderson1986preliminary}, the interaction between individuals with many and few contacts respectively are stochastic, and therefore do not correspond to a fixed (time-invariant) contact graph. For the case of a fixed contact graph, the stability bound~\eqref{eq:stability het} obtained for \eqref{eq:het model} would be conservative because of the assumption that the disease can spread to any contact of an infectious individual. In particular, an infectious individual cannot reinfect its own infector. Therefore, when considering epidemics on configurator model graphs, the shift of the epidemic threshold needs to be corrected by $-1$, resulting in $h=\sigma^2/\mu+\mu-1$, as further discussed in \cite{newman2018,meyers2007contact,kiss2017}.

Relatedly, clustering of the contact network result in the bound \eqref{eq:stability het} becoming conservative. When many neighbors in the network are shared, (e.g. as in small-world networks) local clustering coefficients in the contact network are high and voids the assumption that all neighbors of an infectious node are susceptible. To analytically quantify how local and global clustering effects affect the stability bound of the case isolation scheme, along the lines of \cite{trapman2007analytical}, requires more complex network models. It falls outside the scope of this contribution, but is a worthwhile direction for future work.


\section{Conclusion}
Although obtaining detailed graph models of (the time-varying) human contact graphs relevant for infectious disease transmissions is generally not tractable, control theoretic analyses can provide qualitative, and to some extent quantitative criteria for the feasibility of strategies aimed at halting disease spread. This has been illustrated here by expressing stability conditions for case isolation schemes in homogeneous and heterogeneous populations, as functions of fundamental epidemiological parameters.


\bibliographystyle{IEEETran}
\bibliography{references}
\vspace{1mm}
\end{document}